\documentclass[12pt]{amsart}
\usepackage{graphicx,natbib,amsmath,amssymb}

\theoremstyle{plain}
\newtheorem{thm}{Theorem}[section] \newtheorem{lem}[thm]{Lemma}
\newtheorem{prop}[thm]{Proposition} \newtheorem{cor}[thm]{Corollary}
\theoremstyle{remark} 
\newtheorem*{MT*}{Main theorem} \newtheorem*{MC*}{Main corollary}
\theoremstyle{definition}

\theoremstyle{remark} 
 \newtheorem*{PMT*}{Proof of
the main theorem}

\newcommand{\Z}{{\mathbb Z}}

\newcommand{\x}{{\times}}

\begin{document}

\title{Representations of Degree Three for Semisimple Hopf Algebras} \author{S. Burciu}
%\ead{smburciu@syr.edu}
%\thanks{This paper is part of the author's thesis to be submitted in partial fulfillment of
%the requirements for the Ph.D. degree in Mathematics at Syracuse University.}
\address{Department of Mathematics, Syracuse University, 215 Carnegie Hall, Syracuse, NY 13244}

\begin{abstract} Let $H$ be a cosemisimple Hopf algebra over an algebraically closed field. It is shown that if
$H$ has a simple subcoalgebra of dimension $9$ and has no simple subcoalgebras of even dimension, then $H$
contains either a grouplike element of order $2$ or $3$, or a family of simple subcoalgebras whose dimensions are
the squares of each positive odd integer. In particular, if $H$ is odd dimensional, then its dimension is
divisible by $3$. \end{abstract}

\maketitle

%%% ----------------------------------------------------------------------
%\begin{keyword}Cosemisimple Hopf algebras, Grothendieck group, \newline Kaplansky's conjecture  \MSC 16W30
%\end{keyword}

\section{Introduction}\label{intro}

Let $H$ be a finite dimensional semisimple Hopf algebra over an algebraically closed field $k$. Kaplansky
conjectured that if $H$ has a simple module of dimension $n$, then $n$ divides the dimension of $H$ \cite[Appendix
2]{Kaplansky}.

Several special cases of Kaplansky's conjecture have been proved: Etingof and Gelaki \cite{EG} proved it under the
additional assumption that $H$ is quasi-triangular and $k$ has characteristic $0$ (for another proof see
\cite{SCHN}); Lorenz proved it under the additional assumption that the character of the simple module is central
in $H^*$ \cite{Lo}.

Dually, the Kaplansky conjecture says that if $H$ is a finite dimensional cosemisimple Hopf algebra that contains
a simple subcoalgebra of dimension $n^2$, then $n$ divides the dimension of $H$. This was verified by Nichols and
Richmond in \cite{NR}, for $n=2$. They proved that if $H$ is cosemisimple and contains a simple subcoalgebra of
dimension $4$, then $H$ contains either a Hopf subalgebra of dimension $2$, $12$ or $60$, or a simple subcoalgebra
of dimension $n^2$ for each positive integer $n$. Their approach is based on the study of the Grothendieck group
of $H$.

In this paper, we give a similar treatment for the case $n=3$. Assume that $H$ is cosemisimple and contains a
simple subcoalgebra of dimension $9$. If $H$ has no simple subcoalgebras of even dimension we prove that $H$
contains either a grouplike element of order $2$ or $3$, or a simple subcoalgebra of dimension $n^2$ for each
positive odd integer $n$. In particular, if $H$ is odd dimensional, then its dimension is divisible by $3$. We
remark that if $H$ is bisemisimple, then the assumption that $H$ has odd dimension automatically implies that the
dimension of each simple subcoalgebra is odd (see the corollary on p.~95 of \cite{KSZ}.)

The basic properties of the Grothendieck group of the category of right $H$-comodules are recalled in Section
\ref{S:Grothendieck}. Section \ref{S:Main} contains our main result, namely Theorem \ref{T:Main}, together with
the required lemmas.

We follow the standard notation found in \cite{Montgomery}.  Algebras and coalgebras are defined over $k$;
comultiplication, counit and antipode are denoted by $\Delta$, $\epsilon$ and $S$ respectively; and the category
of right comodules over $H$ is denoted by ${\mathcal M}^H$.

\section{The Grothendieck group ${\mathcal G}(H)$}\label{S:Grothendieck}

Let $H$ be a $k$-coalgebra. For $V \in \mathcal{M}^H$ define $[V]$ to be the isomorphism class of the comodule
$V$. The Grothendieck group of the category of right comodules over $H$ is the free abelian group generated by the
isomorphism classes of simple right $H$-comodules. We denote this group by ${\mathcal G}(H)$. The following result
appears in \cite{NR}.
\begin{prop}\label{descript}
Let $\Gamma$ denote the set of simple subcoalgebras of the coalgebra $H$. For each $C \in \Gamma$, let $V_C$ be a
simple right $C$-comodule. Then ${\mathcal G}(H)$ is the free abelian group with basis $B = \{ [V_C] :\; C \in
\Gamma \}$.
\end{prop}
Here $[V_C]$ is called a basic element and $B$ is called the standard basis. Any basic element $x$ is associated
with a simple subcoalgebra $x_C$ of $H$, and reciprocally, any simple subcoalgebra $C$ of $H$ is associated with a
basic element $C_x$, as in the previous proposition. Every $z \in {\mathcal G}(H)$ may be written uniquely as $z =
\sum_{x \in B} m(x, z) x$ where $m(x, z) \in \Z$.  The integer $m(x,z)$ is called the multiplicity of $x$ in $z$.
If $m(x, z) \ne 0$, then $x$ is called a basic component of $z$. The  multiplicity function may be extended to a
biadditive function $m: {\mathcal G}(H) \x {\mathcal G}(H) \to \Z$ by defining $m(w, z) = \sum_{x \in B} m(x, w)
m(x, z)$. Note that there is a bijection between $\Gamma$ and $B$, the set of isomorphism classes of simple right
$H$-comodules.

If $H$ is a bialgebra, ${\mathcal G}(H)$ becomes a ring. Recall that if $M$ and $N$ are right $H$-comodules then
$M \bigotimes N$ is also a right $H$-comodule via
$$
\rho(m \otimes n) = \sum_{(m),(n)} m_0 \otimes n_0 \otimes m_1 n_1.
$$

Let $1$ denote $[k1]$, the unit of $\mathcal{G}(H)$.  If $V \in \mathcal{M}^H$ is a simple comodule then the
degree of the element $[V] \in \mathcal{G}(H)$ is defined to be the dimension of  the comodule $V$ and it is
denoted by $|[V]|$. Since $B = \{ [V_C] :\; C \in \Gamma \}$ is a basis for $\mathcal{G}(H)$ we get a linear map
called the degree map and defined as above on the canonical basis of $\mathcal{G}(H)$ and then extended by
linearity. By convention, an element $x \in {\mathcal G}(H)$ is said to be $n$-dimensional if its degree is equal
to $n$. Note that $|xy| = |x||y|$ for all $x,y \in B$, and so $|wz|=|w||z|$ for all $w, z \in {\mathcal G}(H)$
which shows that the degree map is a ring homomorphism. Thus $|w||z| = \sum_x m(x , wz)|x|$ for all $w, z \in
{\mathcal G}(H)$. Let $\mathcal{G}(H)_{+}$ be the subset of $\mathcal{G}(H)$ consisting of all elements $z \in
\mathcal{G}(H)$ with the property that $m(x,\; z) \geq 0$ for any basic element $x \in \mathcal{G}(H)$. When $H$
is cosemisimple, the elements of ${\mathcal G}(H)_{+}$ are in a bijective correspondence with all the isomorphism
classes of $H$-comodules $V \in \mathcal{M}^H$ and ${\mathcal G}(H)_{+}$ is closed under multiplication.

Suppose $H$ is a Hopf algebra with antipode $S$.  If $U \in {\mathcal M}^H$ is a right comodule with the comodule
structure $\rho: U \rightarrow U \bigotimes H$ then $U^* \in \mathcal{M}^{H}$ in the following way. Let $\{ u_j
\}$ be a basis of $U$, and write $\rho(u_j) = \sum_i u_i \otimes a_{ij}$, $a_{ij} \in H$. Then if $\{u_i^*\}$ is
the dual basis of $\{u_i\}$ the map $\check{\rho} (u^*_i) = \sum_j u^*_j \otimes S(a_{ij})$ defines a right
$H$-comodule structure on $U^*$. It is easy to see that $C_{U^*} = S(C_U)$ and if $S$ is injective then $U^*$ is
simple whenever $U$ is a simple comodule. Moreover $U^{**} \cong U$ if $S^2(C)=C$. The map $^* : {\mathcal G}(H)
\rightarrow {\mathcal G}(H)$ given by $[M]^*=[M^*]$ is a group homomorphism and a ring antihomomorphism. If $H$ is
cosemisimple then $^*$ is an involution on ${\mathcal G}(H)$. A standard subring in ${\mathcal G}(H)$ is a subring
of ${\mathcal G}(H)$ which is spanned as an abelian group by a subset of $B$.

\begin{thm}\label{T:NR0} (\cite{NR},Theorem 6.) Let $H$ be a bialgebra. There is a one-to-one correspondence
between standard subrings of ${\mathcal G}(H)$ and subbialgebras of $H$ generated as algebras by their simple
subcoalgebras, given by: the subbialgebra $A$ generated by its simple subcoalgebras corresponds to the standard
subring spanned by $\{x_C : C \text{ is a simple subcoalgebra of } A \}$

Furthermore, if $H$ is a cosemisimple Hopf algebra this one-to-one correspondence induces a one-to-one
correspondence between standard subrings closed under `` $^*$'' and Hopf subalgebras of $H$. \end{thm}

The following facts about the multiplicity in ${\mathcal G}(H)$ will be used in the sequel.

\begin{thm}\label{T:NR} (\cite{NR}, Theorem 10.) Let $H$ be a Hopf algebra. \begin{enumerate} \item If the antipode
of $H$ is injective then $m(x,y)=m(x^*, y^*)$ \newline for all $x, y \in {\mathcal G}(H)$. \item If H is
cosemisimple and $k$ is algebraically closed, then \begin{enumerate} \item $m(x,yz)=m(y^*,zx^*)=m(y,xz^*)$ for all
$x, y, z \in {\mathcal G}(H)$ \item For each grouplike element $g$ of $H$, we have $m(g,xy)=1$, if $y=x^* g$ and
$0$ otherwise. \item Let $x \in {\mathcal G}(H)$ be a basic element. Then for any grouplike element $g$ of $H$,
$m(g,xx^{\ast }) > 0$ iff $m(g,xx^{\ast })=1$ iff $gx=x$. The set of such grouplike elements forms a group, of
order at most $|x|^{2}$.
\end{enumerate} \end{enumerate} \end{thm}

For a basic element $x$, the grouplike elements entering in the basic decomposition of $xx^*$ form a subgroup $G$
of the grouplike elements of $H$ (by 2(c) of the previous theorem). Then $kG$ is a Hopf subalgebra of $H$, and by
the freeness theorem \cite{NZ} the order of $G$ divides the dimension of $H$.

\section{The $3$-dimensional case}\label{S:Main}

The goal of this section is to prove the following theorem.

\begin{thm}\label{T:Main}
Let $H$ be a cosemisimple Hopf algebra over an algebraically closed field.  Assume that $H$ contains a simple
subcoalgebra $C$ of dimension $9$ and has no simple subcoalgebras of even dimension. Then one of the following
conditions must hold:

(i) $H$ contains a grouplike element of order $2$ or $3$,

(ii) $H$ has two families of subcoalgebras $\{ C_{2n+1} : n \geq 1 \}$ and $\{ D_{2n+1} : n \geq 1 \}$ with
$\mathrm{dim} C_{2n+1} = \mathrm{dim} D_{2n+1} = (2n+1)^2$ such that
$$
C_{2n+1}C_{3} = C_{2n-1}+D_{2n+1}+C_{2n+3}.
$$
for each $n \geq 1$.
\end{thm}

We need the following lemmas.

\begin{lem}\label{L:basic-decomposition}
Let $k$ be an algebraically closed field and $H$ a cosemisimple Hopf algebra over $k$ with a simple subcoalgebra
$C$ of dimension $9$. Assume that $H$ has no simple subcoalgebras of even dimension. Then one of the following
conditions holds:

\renewcommand{\theenumi}{\roman{enumi}}
\begin{enumerate}
\item $H$ contains a grouplike element of order $2$ or $3$,
\item for any basic element $x_3 \in {\mathcal G}(H)$ with $|x_3|=3$ and $x_3 x_3^*=1+u+v$
where $u$ and $v$ are basic elements of ${\mathcal G}(H)$ with $|u|=3$ and $|v|=5$.
\end{enumerate}
\end{lem}

\begin{proof}
Suppose that $x_3$ is a basic element of degree $3$. Then $x_3 x_3^* = 1+y$ where $|y|=8$. Note that since $H$ has
no simple subcoalgebras of even dimension, $y$ is not a basic element. We consider the decomposition of $y$ into
basic elements in ${\mathcal G}(H)$. If $y$ has $1$, $2$, $3$, $5$ or $8$  one dimensional representations in this
decomposition, then together with $1$ they form a subgroup of ${\mathcal G}(H)$ with $2$, $3$, $4$, $6$ or $9$
elements respectively, and so (i) holds.

Since $|y|=8$ and there are no $2$-dimensional simple comodules of $H$, it is clear that $y$ cannot have exactly
$6$ grouplike elements. A similar argument shows that $y$ cannot have exactly $4$ grouplike elements.

It follows that if (i) does not hold then $y$ has no grouplike elements in its basic decomposition.  Since any
other element of ${\mathcal G}(H)$ has degree at least $3$, $y=u+v$, where $u$ and $v$ are basic elements with
$|u|=3$ and $|v|=5$. In this case (ii) holds.
\end{proof}

\begin{lem}\label{L:selfdual}
If $H$ satisfies the assumptions of Theorem \ref{T:Main}, then one of the following conditions holds:
\renewcommand{\theenumi}{\roman{enumi}}
\begin{enumerate}
\item $H$ contains a grouplike element of order $2$ or $3$,
\item there is a $3$-dimensional self dual basic element $x_3$
with $x^2_3 = 1+x_3+x_5$, where $x_5$ is a $5$-dimensional basic element.
\end{enumerate}
\end{lem}

\begin{proof}
Assume (i) does not hold. By the previous lemma there is a $3$-dimensional basic element with $x_3 x_3^* = 1+u+v$,
where $u$ and $v$ are basic elements of ${\mathcal G}(H)$ with $|u|=3$ and $|v|=5$. Since $x_3 x_3^*$ is self
adjoint, the last relation implies that $u=u^*$ and $v=v^*$. If $u=x_3$, we are done. Otherwise the previous lemma
applied to $u$ instead of $x_3$ gives that $u u^* = 1+u_1+v_1$ with $|u_1|=3$, $|v_1|=5$, and $u_1$, $v_1$
self-adjoint. If $u \neq u_1$, Lemma \ref{L:basic-decomposition} again gives $u_1 u_1^* = 1+u_2+v_2$ where $u_2$,
$v_2$ are self-adjoint with $|u_2|=3$ and $|v_2|=5$. It suffices to show that $u_1=u_2$.

Since $m(u_1, u^2)=1$, it follows from Theorem \ref{T:NR} that $m(u, u u_1)=1$. Suppose $u u_1=u+y$. Since
$|y|=6$, $y$ is not a basic element of ${\mathcal G}(H)$. Therefore $y=w+\xi$ where $w,\xi \neq 0$.

Using the previous relations we have
$$
u^2 u_1=(1+u_1+v_1) u_1=u_1+1+u_2+v_2+v_1 u_1.
$$
On the other hand
$$
u^2 u_1=u (u+w+\xi)=1+u_1+v_1+u w+u \xi.
$$
Hence
\begin{equation}\label{eqn1}
u_2+v_2+v_1 u_1=v_1+u w+u \xi.
\end{equation}
But $1\leq m(w, u u_1)=m(u_1, u w)$. Assume $u_1 \neq u_2$. In this case, from the last two relations, $u_1$ must
enter in the basic decomposition of $v_1 u_1$ and $1 \leq m(u_1, v_1 u_1)=m(v_1, u^2_1)$. Therefore $v_1=v_2$ and
$m(u_1,v_1 u_1)=m(v_1, u_1^2)=1$. Then (\ref{eqn1}) becomes $v_1 u_1+u_2=u w+u \xi$. This relation is impossible
in the case $u_1 \neq u_2$.  Indeed, if $u_1 \neq u_2$ then the multiplicity of $u_1$ on the left hand side is
equal to $1$ whereas on the other side the multiplicity of $u_1$ is at least $2$ since $u_1$ enters in both terms
of the sum.
\end{proof}

Before giving the proof of Theorem \ref{T:Main}, one more lemma is needed.

\begin{lem}\label{L:triple}
Suppose $H$ satisfies the assumptions of Theorem \ref{T:Main} and \newline $a, a', b \in {\mathcal G}(H)$ with
$|a|=|a'|<|b|$ and $a, b$ basic elements. Let $x_3$ be a $3$-dimensional basic element of ${\mathcal G}(H)$ with
$a x_3=b+c+u=a'+c+v$, for some $c, u, v \in {\mathcal G}(H)_{+}$ with $c \neq 0$. Then $v=b$ and $u=a'$.
\end{lem}

\begin{proof} It is easy to see that any basic component of $a x_3$ has degree at least $|a|/3$. Indeed, if $1 \leq
m(z, a x_3) = m(a, z x_3)$ then $|a| \leq |z x_3|$ and $|z| \geq |a|/3$. We have $b+u=a'+v$. Since $|a'|<|b|$ and
$b$ is a basic element, it follows that $b$ is a basic component of $v$. Therefore $v=v_1 + b$ and $u=a'+v_1$. We
want to show that $v_1=0$. If $v_1 \neq 0$ then $c+v_1$ must have at least three basic components since every
basic element has odd degree. These are also basic components of $a x_3$.  But this is impossible, since in that
case $|c+v_1| \geq 3|a|/3=|a|$ and therefore the degree of $b+c+u$ is strictly greater than $3|a|$. Thus $v_1= 0$.
\end{proof}

We are ready to prove our main result.
\begin{proof} (of Theorem \ref{T:Main}.)
Assume (i) does not hold. By Lemma \ref{L:selfdual} there is a $3$-dimensional basic element $x_3$ such that
$x_3^2=1+x_3+x_5$.

It will be shown that (ii) holds in this case. For (ii) it suffices to prove the existence of two families of
basic elements $\{x_{2n+1}:n\geq 1\}$, $\{x'_{2n+1}: n\geq 1\}$ corresponding to the two families of simple
subcoalgebras $\{C_{2n+1}:n\geq 1\}$, $\{D_{2n+1} : n\geq 1\}$ satisfying $|x_{2n+1}|=|x'_{2n+1}|=2n+1$ and
$$
x_{2n+1}x_3=x_{2n-1}+x'_{2n+1}+x_{2n+3}
$$
for all $n\geq 0$.

For $x_1= 1$ and $x'_3=x_3$ the first relation of (ii) is satisfied. Suppose we have found $x_3, x_5, ...
x_{2n+1}, x_{2n+3}$ and $ x'_3, x'_5, ... x'_{2n+1}$ such that:
$$
x_{2k+1}x_3=x_{2k-1}+x'_{2k+1}+x_{2k+3}
$$
for any natural number $k$ with $1\leq k\leq n$.

We want to show that there are another two basic elements $x'_{2n+3}$, $x_{2n+5}$ such that
$$
x_{2n+3}x_3=x_{2n+1}+x'_{2n+3}+x_{2n+5}
$$
Since $m(x_{2n+1},x_{2n+3}x_3)=m(x_{2n+3},x_{2n+1}x_3)=1$ we may write
$$
x_{2n+3}x_3=x_{2n+1}+y_0+z_0,
$$
where $y_0$ is a basic component of $x_{2n+3}x_3$, different from $x_{2n+1}$ and with the smallest possible
degree. If $|y_0|\geq 2n+3$ we are done. Indeed any other irreducible that enters in $z_0$ has dimension at least
$2n+3$. If $z_0$ is not basic then it contains at least three basic elements since $|z_0|$ is odd. In this case
$|z_0|\geq 3(2n+3)$ which is not possible. Hence $z_0$ is basic. Since $|y_0| + |z_0|= 4n+8$ and $|z_0|\geq
|y_0|\geq 2n+3$ the following equalities are satisfied $|y_0|= 2n+3$ and $|z_0|=2n+5$. Then let $x'_{2n+3}=y_0$
and $x_{2n+5}=z_0$.

The case $|y_0| < 2n+3$ will be shown to be impossible. Note that $m(x_{2n+3}, y_0x_3)=m(y_0, x_{2n+3}x_3)=1$ and
if $y_0 \neq x_{2n+3}x_3$ we may write again
$$
y_0x_3=x_{2n+3}+y_1+z_1,
$$
where $y_1$ is a basic component of $y_0x_3$ different from $x_{2n+3}$ and with the smallest possible dimension.
The degree of $y_0x_3-x_{2n+3}$ is even and therefore there are at least three basic components in the
decomposition of $y_0x_3$. Since $|y_0|<2n+3$ we have $|y_1|<|y_0|$. The same procedure gives $y_1
x_3=y_0+y_2+z_2$ where $y_2$ is again a basic component of $y_1x_3$, different from $y_0$ and with the smallest
possible dimension. Similarly $|y_2|<|y_1|$.

In this manner we construct a sequence of basic elements $y_0$, $y_1, ... y_k$ with
\newline $|y_k|<|y_{k-1}|<...<|y_2|<|y_1|<|y_0|<2n+3$ and
\begin{eqnarray*}
y_0x_3=&x_{2n+3}+y_1+z_1   \\
y_1x_3=&y_0+y_2+z_2   \\
\vdots & \\
y_{k-1}x_3=&y_{k-2}+y_k+z_k
\end{eqnarray*}
where the $z_i$ are not necessarily irreducible.

Since the dimension of $y_k$ is decreasing this process must stop. Therefore we may suppose $y_k x_3$ is basic, so
$y_k x_3= y_{k-1}$. Note that the case $y_0x_3=x_{2n+3}$ simply means that the process stops after the first
stage.

It will be shown that $k=n=1$. First note that $k<n+1$ since $|y_k|<|y_{k-1}|<...<|y_1|<|y_0|<2n+3$ and all the
elements have odd degree. Next we will prove that
\begin{equation} \label{eqn2} y_{k-t}= y_k x_{2t+1} \end{equation}
for $1\leq t\leq k+1$. For consistency of notation we put $y_{-1}=x_{2n+3}$.

For $t=1$, $y_k x_3=y_{k-1}$ from above. Suppose $$ y_{k-t} =  y_k x_{2t+1} $$ for $1\leq t\leq s$. We need to
prove that $$ y_{k-s-1}= y_kx_{2s+3} $$ For $t=s$ we have $y_{k-s}= y_k x_{2s+1}$. Multiplication by $x_3$ gives
$y_{k-s} x_3= y_k x_{2s+1}x_3$. Then $$ y_{k-s-1}+y_{k-s+1}+z_{k-s+1}= y_k x_{2s-1}+y_k x'_{2s+1}+y_k x_{2s+3} $$
(We used that $2s+1\leq 2k+1 \leq 2n+1$). Since $y_{k-s+1} = y_k x_{2s-1}$ it follows that
$$ y_{k-s-1}+z_{k-s+1}= y_k x'_{2s+1}+y_kx_{2s+3} $$ But $|y_k
x'_{2s+1}|=|y_k x_{2s+1}|=|y_{k-s}|<|y_{k-s-1}|$. Lemma \ref{L:triple} applied to  $a=y_{k-s}$, $a'=y_k x'_{2s+1}$
and $b=y_{k-s-1}$ implies that $z_{k-s+1}=y_kx'_{2s+1}$ and $y_{k-(s+1)}= y_kx_{2s+3}$ as required. Note that
$y_{-2}$ represents $x_{2n+1}$ in the case $t=k+1$.

For $t=k+1$ relation (\ref{eqn2}) becomes
\begin{equation}\label{eqn3} x_{2n+3}=y_k x_{2k+3} \end{equation}
We show that this is impossible if $k<n$. Indeed, $2k+3<2n+3$ and
\newline relation (\ref{eqn3}) multiplied by $x_3$ gives \begin{equation}\nonumber
x_{2n+1}+y_0+z_0 = y_kx_{2k+1}+y_k x'_{2k+3}+y_k x_{2k+5}\end{equation}  Since $y_0=y_kx_{2k+1}$ it follows that
\begin{equation}\label{eqn4} x_{2n+1}+z_0= y_k x'_{2k+3}+y_k x_{2k+5}\end{equation} It will be shown that $x_{2n+1}$ cannot be a
basic component of either term of the right hand side. The degree of the first term is  $$ |y_k x'_{2k+3}|=|y_k
x_{2k+3}|=|x_{2n+3}|=2n+3$$ and the difference between its degree and the degree of $x_{2n+1}$ is $2$. The other
components of $y_k x'_{2k+3}$ cannot be $1$-dimensional since they are also components of $x_{2n+3} x_3$. If
$x_{2n+1}$ is a basic component of the second term, $y_k x_{2k+5}$, then relation \ref{eqn4} implies that $y_k
x_{2k+5}=x_{2n+1}+u+v$, $z_0=y_k x'_{2k+3}+u+v$ and $x_{2n+3} x_3=x_{2n+1}+y_0+y_k x'_{2k+3}+u+v$. Therefore
$|y_0|+|u|+|v|=2n+5$. Since any basic components of $x_{2n+3}x_3$ has degree at least $\frac{2n+3}{3}$ and $y_0$
is the smallest component different from $x_{2n+1}$, we have $\frac{2n+3}{3} \leq |y_0| \leq \frac{2n+5}{3}$. But
$y_0x_3$ contains $x_{2n+3}$ and $|y_0x_3|-|x_{2n+3}|\leq 2$. So if $y_0 x_3-x_{2n+3}\neq 0$ then it has to be the
sum of two grouplike elements, which is impossible. Indeed, if $g$ and $h$ were these two grouplike elements then
$gy_0=hy_0=x_3$ and $gh^{-1}x_3=x_3$ which implies $g=h$ and contradicts part 2 (c) of Theorem~\ref{T:NR}. We
conclude that $y_0 x_3=x_{2n+3}$, and hence $|y_0|= \frac{2n+3}{3}$. Moreover, $|u|=|y_0|=\frac{2n+3}{3}$ and
$|v|=|y_0|+2$. Thus $u x_3=x_{2n+3}$ and $v x_3=x_{2n+3}+w$, where $|w|=6$. But $n\geq 1$, $|v| \geq 4$ and so $v
x_3$ cannot have any $1$-dimensional basic components. Therefore $w$ is the sum of two $3$-dimensional basic
components $w_1$ and $w_2$. Since $1 \leq m(w_i, vx_3)=m(v,w_ix_3)$ and $|w_ix_3|=9$,  it follows that $w_i x_3 -
v$ has dimension at most $5$. Multiplying on the right by $x_3$, it is easy to see that $w_i x_3-v \neq 0$. So
each product $w_i x_3$ has a $1$-dimensional component $g_i$ with $w_i=g_ix_3$. Therefore $v x^2_3=x_{2n+3}
x_3+w_1 x_3+w_2 x_3
        =x_{2n+3} x_3+g_1 x^2_3+g_2 x^2_3$
has two components of degree $1$, namely $g_1$ and $g_2$. On the other hand, $v x_3^2=v+v x_3+v x_5$ and the only
$1$-dimensional components of $v x_3^2$ appear in the basic decomposition of $v x_5$. Part 2(b) of Theorem
\ref{T:NR} implies that $v=g_1 x_5=g_2 x_5$. So for $i=1,2$, $v x_3=g_i x_5 x_3$ and its basic components have
dimension at most $7$.  Since $2n+3$ is divisible by $3$ and $n\geq 1$ we have $|x_{2n+3}|=2n+3\geq 9$. Hence
$x_{2n+3}$ cannot be a basic component of $v x_3$. This means that the relation (\ref{eqn3}) cannot hold for
$k<n$.

Now suppose $k=n$. Then $x_{2n+3}=y_n x_{2n+3}$ implies that $|y_{n}|=1$ and $y_n=g$, a grouplike element of $H$.
Moreover, relation \ref{eqn2} implies that $y_{n-t}=gx_{2t+1}$, for $0\leq t\leq n$. Then
\begin{equation}\label{eqn5}
x_{2n+3}x_3=x_{2n+1}+gx_{2n+1}+z_0,
\end{equation}
which gives $|z_0|=2n+7$.  We have to consider two cases, whether $z_0$ is a basic element or not.

If $z_0$ is basic then $g$ has order $2$, which would imply (i), contrary to our assumption. Indeed, the last
relation multiplied by $g$ becomes $x_{2n+3}x_3=gx_{2n+3}x_3=gx_{2n+1}+g^2x_{2n+1}+gz_0$. Therefore
$g^2x_{2n+1}=x_{2n+1}$. The decomposition formula for $x_{2n+1}x_3$ multiplied by $x_3$ gives
$g^2x_{2n-1}=x_{2n-1}$. Similarly, $g^2x_{2t+1}=x_{2t+1}$ for any $0\leq t\leq n$. In particular $g^2 x_1=x_1$
gives $g^2=1$.

Hence $z_0$ cannot be a basic element. In this case, in its basic decomposition there are at least three terms,
since every basic element has odd degree. By the choice of $y_0$, any of these basic elements has degree at least
$|y_0|=2n+1$; therefore $2n+7\geq 3(2n+1)$ which implies $n=1$.

It will be shown that this is impossible. Write $z_0 = u_1+u_2+u_3$.  Then (\ref{eqn5}) becomes $$ x_5 x_3 = x_3 +
g x_3 + u_1 + u_2 + u_3. $$ and one dimensional basic elements cannot appear on the right hand side. Hence
$|u_i|=3$ for $1\leq i\leq 3$. It is easy to see that each of the products $u_i x_3$ has a one dimensional basic
component since each of them has degree $9$ and all of them have $x_5$ as a component. Hence, by part 2 (b) of
Theorem \ref{T:NR}, $u_i=h_i x_3$ where each $h_i$ is a grouplike element.

Consider the set $V$ of grouplike elements $h$ such that $h x_5=x_5$. By part 2 (c) of Theorem \ref{T:NR}, $V$ is
a group. It is easy to see that $1, g, h_1, h_2, h_3 \in V$. If $h\in V$ then $m(hx_3, x_5x_3) = m(hx_5, x^2_3) =
m(x_5, x^2_3)=1$, implying that $hx_3$ is one of the basic elements $u_i$. But if $hx_3=h_ix_3$ then
multiplication by $x_3$ on the right gives $h+hx_3+hx_5=h_i+h_ix_3+h_ix_5$ and so $h=h_i$. Therefore $V=\{ 1, g,
h_1, h_2, h_3 \}$. The relation $x_3^2=1+x_3+x_5$ multiplied by $x_5$ on the right gives $x_5^2 =
4x_5+1+g+h_1+h_2+h_3$. This shows that $\{x_3, x_5\} \cup  V$ generates a standard subring $R_1$ of ${\mathcal
G}(H)$. We have $x_3^3=x_3^2 x_3=x_3+x_3^2+x_5 x_3$. On the other hand, $x_3^3=x_3 x_3^2=x_3+x_3^2+x_3 x_5$.
Therefore $x_5 x_3=x_3 x_5$. But $x_5 x_3=(x_3 x_5)^*$ and $R_1$ is closed under $^*$. By Theorem \ref{T:NR0}, it
corresponds to a Hopf subalgebra $K_1$ of $H$ with dimension equal to $5*1^2 + 5*3^2 + 5^2 =75$. In a similar way
it can be seen that $\{x_5\} \cup V$ generates a standard subring $R_2$. It is also closed under $^*$ and by the
same argument it corresponds to a Hopf subalgebra $K_2$ of $H$ with dimension $30$. Clearly $R_2 \subseteq R_1$.
Then $K_2$ is a Hopf subalgebra of $K_1$ contradicting the freeness theorem.
\end{proof}

\begin{cor}
Let $H$ be an odd dimensional Hopf algebra over an algebraically closed field $k$. If $H$ contains a simple
subcoalgebra of dimension $9$ and no simple subcoalgebras of even dimension then $\mathrm{dim}_{k} H$ is divisible
by $3$.
\end{cor}

\begin{proof}
Since $\mathrm{dim}_k H$ is odd the freeness theorem implies that $H$ cannot have a grouplike element of order
$2$. Therefore $H$ has a grouplike element of order $3$ and the same theorem gives the divisibility relation.
\end{proof}

A similar result was obtained by Kashina, Sommerh\"auser and Zhu in \cite{KSZ'}. They assume that the
characteristic of the base field is $0$ in which case the assumption that $H$ is odd dimensional automatically
implies that $H$ has no even dimensional simple subcoalgebras.

\subsection*{Acknowledgment} I would like to express my deep gratitude
to my thesis supervisor Declan Quinn whose guidance and support were crucial for the successful completion of this
paper. I also thank Y. Sommerh\"auser for his encouraging comments and suggestions regarding the presentation of
this paper.

\bibliographystyle{plain} \bibliography{ggha3matharchive}
\end{document}